# EMPIRICAL-LIKELIHOOD-BASED CONFIDENCE INTERVAL FOR THE MEAN WITH A HEAVY-TAILED DISTRIBUTION

By Liang Peng

*Georgia Institute of Technology*

Empirical-likelihood-based confidence intervals for a mean were introduced by Owen [*Biometrika* **75** (1988) 237–249], where at least a finite second moment is required. This excludes some important distributions, for example, those in the domain of attraction of a stable law with index between 1 and 2. In this article we use a method similar to Qin and Wong [*Scand. J. Statist.* **23** (1996) 209–219] to derive an empirical-likelihood-based confidence interval for the mean when the underlying distribution has heavy tails. Our method can easily be extended to obtain a confidence interval for any order of moment of a heavy-tailed distribution.

**1. Introduction.** Suppose $X_1,\ldots,X_n$ are i.i.d. random variables with common distribution function $F$ satisfying

$$
\begin{aligned}
1 - F(x) &= x^{-\alpha_R} L_1(x), \\
F(-x) &= x^{-\alpha_L} L_2(x),
\end{aligned}
\tag{1.1}
$$

where $\alpha_R > 1$, $\alpha_L > 1$ and $L_1(x)$ and $L_2(x)$ are slowly varying functions, that is, $\lim_{t\to\infty} L_1(tx)/L_1(t) = 1$ and $\lim_{t\to\infty} L_2(tx)/L_2(t) = 1$ for all $x > 0$. We are interested in obtaining a confidence interval for the mean $\mu = E(X_i)$. In the case $\alpha_R = \alpha_L = \alpha \in (1,2)$ and $\lim_{x\to\infty} (L_1(x))/(L_1(x) + L_2(x)) = p \in [0,1]$, we know that $F$ lies in the domain of attraction of a stable law with index $\alpha$. If both $\alpha_R$ and $\alpha_L$ are greater than 2, then $F$ is in the domain of attraction of a normal distribution. The use of heavy-tailed models in financial markets, such as value-at-risk in the context of risk management, is attracting more attention; see the book by Embrechts, Klüppelberg and Mikosch (1997). There are many other applications of heavy-tailed distributions, for example, teletraffic data [see Resnick (1997)] and community size [see Feuerverger and Hall (1999)].

To estimate the population mean $\mu$, a natural estimator is the sample mean $\bar{X}_n = \frac{1}{n}\sum_{i=1}^n X_i$. Since the type of limiting distribution of the sample mean, a stable law or a normal, depends on the tail indices $\alpha_L$ and $\alpha_R$, it is hard to derive a confidence interval for $\mu$ based on the limit, especially when the limit is a stable law with index between 1 and 2. To avoid distinguishing the type of limit, Hall and LePage (1996) proposed a subsample bootstrap method to construct a confidence interval for $\mu$ under more general conditions than (1.1). However, this method performed poorly in the case $\alpha_L = \alpha_R \in (1,2)$ [see Hall and Jing (1998)]. More recently, Peng (2001) proposed a new estimator for $\mu$ that has a limit that is always a normal distribution by assuming $\alpha_R = \alpha_L > 1$ and second order regular variation conditions for $1 - F(x)$ and $F(-x)$. This result prompted us to look at the empirical likelihood method. Empirical likelihood methods for constructing confidence regions were introduced by Owen (1988, 1990). One of their advantages is that they enable the shape of a region, such as the degree of asymmetry in case of a confidence interval, to be determined automatically by the sample. In certain regular cases, empirical-likelihood-based confidence regions are Bartlett correctable; see Hall and La Scala (1990), DiCiccio, Hall and Romano (1991) and Chen and Hall (1993). For more details on empirical likelihood methods, see Owen (2001).

We organize this article as follows. In Section 2 we introduce our empirical likelihood method and state our main results. A simulation study and a real application are given in Section 3. All proofs are deferred to Section 4.

**2. Methodology and main results.** The key idea of our method is to employ extreme value theory to deal with two tails of the underlying distribution and to treat the middle part of the underlying distribution nonparametrically. Since this approach involves estimation of tail index, we know we need a stricter condition than (1.1) to investigate the asymptotic behavior of the tail index estimator. For simplicity we assume, as $x \to \infty$,

$$
\begin{aligned}
1 - F(x) &= c_R x^{-\alpha_R}\{1 + b_R x^{-\beta_R} + o(x^{-\beta_R})\}, \\
F(-x) &= c_L x^{-\alpha_L}\{1 + b_L x^{-\beta_L} + o(x^{-\beta_L})\},
\end{aligned}
\tag{2.1}
$$

where $c_R > 0$, $c_L > 0$, $b_R \neq 0$, $b_L \neq 0$, $\beta_R > 0$ and $\beta_L > 0$. We remark that (2.1) is a special case of a general second-order regular variation [see de Haan and Stadtmüller (1996)].

Put $\delta_i = I(X_i \leq X_{n,m})$ and $\bar{\delta}_i = I(X_i \geq X_{n,n-k+1})$, where $X_{n,1} \leq \cdots \leq X_{n,n}$ denote the order statistics of $X_1, \ldots, X_n$, $m = m(n) \to \infty$, $m/n \to 0$, $k = k(n) \to \infty$, $k/n \to 0$, $X_{n,m} \leq 0$ and $X_{n,n-k+1} \geq 0$. We may approximate the two tails of $F$ by

$$F(x) = 1 - c_R x^{-\alpha_R}, \qquad \text{if } x \geq X_{n,n-k+1},$$

and
$$F(x) = c_L(-x)^{-\alpha_L}, \qquad \text{if } x \leq X_{n,m}.$$

Thus the log-likelihood function may be written as two parts, the parametric log-likelihood based on data in two tails and the nonparametric likelihood based on data in the midrange, that is,
$$l_0(c_R, c_L, \alpha_R, \alpha_L, p_i) = l_1(c_R, c_L, \alpha_R, \alpha_L) + l_2(p_i),$$
where
$$\begin{aligned} l_1(c_R, c_L, \alpha_R, \alpha_L) \\ = k \log \alpha_R + k \log c_R - (\alpha_R + 1) \sum_{\bar{\delta}_i=1} \log X_i \\ + m \log \alpha_L + m \log c_L - (\alpha_L + 1) \sum_{\delta_i=1} \log(-X_i) \end{aligned}$$
and
$$l_2(p_i) = \sum_{\bar{\delta}_i=\delta_i=0} \log p_i.$$

By considering the log-likelihood function $l_0(c_R, c_L, \alpha_R, \alpha_L, p_i)$ subject to the constraints

(2.2) $\qquad \alpha_R > 1, \qquad \alpha_L > 1, \qquad c_R > 0, \qquad c_L > 0, \qquad p_i > 0,$

(2.3) $\qquad \sum_{\bar{\delta}_i=\delta_i=0} p_i = 1 - c_R X_{n,n-k+1}^{-\alpha_R} - c_L(-X_{n,m})^{-\alpha_L},$

(2.4) $\qquad \sum_{\bar{\delta}_i=\delta_i=0} p_i X_i = \mu - c_R \frac{\alpha_R}{\alpha_R - 1} X_{n,n-k+1}^{1-\alpha_R} - c_L \frac{\alpha_L}{\alpha_L - 1}(-X_{n,m})^{1-\alpha_L},$

we can employ the semiparametric likelihood ratio method, as in Qin and Wong (1996), to obtain a confidence interval for $\mu$. The details are given in the following discussion.

First we maximize $l_0$ subject to constraints (2.2) and (2.3). By the method of Lagrange multipliers, we easily obtain

(2.5)
$$\bar{\alpha}_R = \left\{ \frac{1}{k} \sum_{\bar{\delta}_i=1} \log \frac{X_i}{X_{n,n-k+1}} \right\}^{-1},$$
$$\bar{\alpha}_L = \left\{ \frac{1}{m} \sum_{\delta_i=1} \log \frac{X_i}{X_{n,m}} \right\}^{-1},$$
$$\bar{c}_R = \frac{k}{n} X_{n,n-k+1}^{\bar{\alpha}_R},$$

$$\bar{c}_L = \frac{m}{n}(-X_{n,m})^{\bar{\alpha}_L},$$

$$\bar{p}_i = n^{-1}, \qquad \text{if } \bar{\delta}_i = \delta_i = 0.$$

Note that both $\bar{\alpha}_R$ and $\bar{\alpha}_L$ are Hill estimators [see Hill (1975)].

Next we maximize $l_0$ subject to constraints (2.2)–(2.4) for fixed $\mu$. Define

$$w_i = w_i(\mu, c_R, c_L, \alpha_R, \alpha_L)$$
$$= X_i - (\mu - c_R \alpha_R (\alpha_R - 1)^{-1} X_{n,n-k+1}^{1-\alpha_R}$$
$$\qquad - c_L \alpha_L (\alpha_L - 1)^{-1} (-X_{n,m})^{1-\alpha_L})$$
$$\qquad \times (1 - c_R X_{n,n-k+1}^{-\alpha_R} - c_L(-X_{n,m})^{-\alpha_L})^{-1},$$

$$\hat{c}_R(\alpha_R, \mu) = \frac{k}{n} X_{n,n-k+1}^{\alpha_R}$$
$$\qquad - \left\{ \sum_{\bar{\delta}_i=1} \log \frac{X_i}{X_{n,n-k+1}} - \frac{k}{\alpha_R} \right\}$$
$$\qquad \times \{\alpha_R(\alpha_R - 1)n^{-1} X_{n,n-k+1}^{\alpha_R}$$
$$\qquad - \mu(\alpha_R - 1)^2 n^{-1} X_{n,n-k+1}^{\alpha_R - 1}\},$$

$$\hat{c}_L(\alpha_L, \mu) = \frac{m}{n}(-X_{n,m})^{\alpha_L}$$
$$\qquad - \left\{ \sum_{\delta_i=1} \log \frac{X_i}{X_{n,m}} - \frac{m}{\alpha_L} \right\}$$
$$\qquad \times \{\alpha_L(\alpha_L - 1)n^{-1}(-X_{n,m})^{\alpha_L}$$
$$\qquad - \mu(\alpha_L - 1)^2 n^{-1}(-X_{n,m})^{\alpha_L - 1}\},$$

$$\hat{\lambda}_2(\alpha_R, \mu) = \left\{ kn^{-1}(\alpha_R - 1)^{-2} X_{n,n-k+1} \left[ \sum_{\bar{\delta}_i=1} \log \frac{X_i}{X_{n,n-k+1}} - \frac{k}{\alpha_R} \right]^{-1} \right.$$
$$\qquad \left. - n^{-1}\alpha_R(\alpha_R - 1)^{-1} X_{n,n-k+1} + n^{-1}\mu \right\}^{-1},$$

$$\bar{\lambda}_2(\alpha_L, \mu) = \left\{ mn^{-1}(\alpha_L - 1)^{-2}(-X_{n,m}) \left[ \sum_{\delta_i=1} \log \frac{X_i}{X_{n,m}} - \frac{m}{\alpha_L} \right]^{-1} \right.$$
$$\qquad \left. - n^{-1}\alpha_L(\alpha_L - 1)^{-1}(-X_{n,m}) + n^{-1}\mu \right\}^{-1},$$

$$\hat{\lambda}_1(\alpha_R, \alpha_L, \mu) = (n - m - k)$$

$$\times \{1 - \hat{c}_R(\alpha_R,\mu)X_{n,n-k+1}^{-\alpha_R} - \hat{c}_L(\alpha_L,\mu)(-X_{n,m})^{-\alpha_L}\}^{-1}.$$

Let $\hat{\alpha}_R = \hat{\alpha}_R(\mu)$ and $\hat{\alpha}_L = \hat{\alpha}_L(\mu)$ denote the solutions to the equations

$$\begin{aligned}
\hat{h}(\alpha_R,\alpha_L,\mu) &= \sum_{\bar{\delta}_i=\delta_i=0} (w_i(\mu,\hat{c}_R(\alpha_R,\mu),\hat{c}_L(\alpha_L,\mu),\alpha_R,\alpha_L)) \\
&\quad \times (\hat{\lambda}_1(\alpha_R,\alpha_L,\mu) \\
&\qquad + \hat{\lambda}_2(\alpha_R,\mu)w_i(\mu,\hat{c}_R(\alpha_R,\mu),\hat{c}_L(\alpha_L,\mu),\alpha_R,\alpha_L))^{-1} \\
&= 0
\end{aligned}$$ (2.6)

and

$$\begin{aligned}
\bar{h}(\alpha_R,\alpha_L,\mu) &= \sum_{\bar{\delta}_i=\delta_i=0} (w_i(\mu,\hat{c}_R(\alpha_R,\mu),\hat{c}_L(\alpha_L,\mu),\alpha_R,\alpha_L)) \\
&\quad \times (\hat{\lambda}_1(\alpha_R,\alpha_L,\mu) \\
&\qquad + \bar{\lambda}_2(\alpha_L,\mu)w_i(\mu,\hat{c}_R(\alpha_R,\mu),\hat{c}_L(\alpha_L,\mu),\alpha_R,\alpha_L))^{-1} \\
&= 0.
\end{aligned}$$ (2.7)

(The existence is given in Lemma 3 in Section 4.) Hence

$$\begin{aligned}
\hat{c}_R &= \hat{c}_R(\hat{\alpha}_R,\mu), \\
\hat{c}_L &= \hat{c}_L(\hat{\alpha}_L,\mu), \\
\hat{\lambda}_1 &= \hat{\lambda}_1(\hat{\alpha}_R,\hat{\alpha}_L,\mu), \\
\hat{\lambda}_2 &= \hat{\lambda}_2(\hat{\alpha}_R,\mu), \\
\hat{p}_i &= \{\hat{\lambda}_1 + \hat{\lambda}_2 w_i(\mu,\hat{c}_R,\hat{c}_L,\hat{\alpha}_R,\hat{\alpha}_L)\}^{-1}, \qquad \text{if } \bar{\delta}_i = \delta_i = 0
\end{aligned}$$ (2.8)

are the values which maximize $l_0$ subject to (2.2)–(2.4). (The proof is given in Section 4.) Therefore, the semiparametric likelihood ratio multiplied by $-2$ is defined as

$$l(\mu) = -2\{l_0(\hat{c}_R,\hat{c}_L,\hat{\alpha}_R,\hat{\alpha}_L,\hat{p}_i) - l_0(\bar{c}_R,\bar{c}_L,\bar{\alpha}_R,\bar{c}_L,\bar{p}_i)\}.$$

Our main results are as follows.

THEOREM 1. *Suppose (2.1) holds. Choose $k = o(n^{2\beta_R/(\alpha_R+2\beta_R)})$, $k/\log n \to \infty$, $m = o(n^{2\beta_L/(\alpha_L+2\beta_L)})$ and $m/\log n \to \infty$ as $n \to \infty$. Further assume $B_0 n^{l_0} \leq m, k \leq B_1 n^{l_1}$ for some $B_0 > 0$, $B_1 > 0$, $0 < l_0 \leq l_1 < 1$ and*

$$\lim_{n\to\infty} \frac{(k/n)^{1-2/\alpha_R}}{(k/n)^{1-2/\alpha_R} + (m/n)^{1-2/\alpha_L}}$$

*exists. Let the true mean be $\mu_0$. Then*

$$l(\mu_0) \xrightarrow{d} \chi^2_{(1)}.$$

Based on Theorem 1, a simple approach to construct a $1 - \alpha$ level confidence interval for $\mu_0$ is

$$I_{1-\alpha} = \{\mu : l(\mu) \leq d_{1-\alpha}\},$$

where $d_{1-\alpha}$ is the $1 - \alpha$ quantile of a $\chi^2_{(1)}$ distribution. This gives a confidence interval for $\mu_0$ with asymptotically correct coverage probability $1 - \alpha$, as stated in the following corollary.

COROLLARY 1. *Assume the conditions in Theorem 1 hold. Then as $n \to \infty$,*

$$P(\mu_0 \in I_{1-\alpha}) = 1 - \alpha + o(1).$$

REMARK 1. Our method can easily be extended to obtain a confidence interval for any order of moment of a heavy-tailed distribution.

REMARK 2. The choices of sample fractions $k$ and $m$ are difficult both theoretically and practically. This will be a part of our future research plan. However, our simulation study in Section 3 shows that this approach is robust against the choices of sample fraction.

REMARK 3. An empirical likelihood estimator for $\mu$ can be obtained as $\hat{\mu} = \arg\min_\mu l(\mu)$, which can be shown to give $\hat{\mu} - \mu = O_p(\frac{\sigma(m/n,k/n)}{\sqrt{n}})$, where $\sigma(s,t)$ is defined in Lemma 1. That is, the convergence rate of $\hat{\mu}$ is the same as that of the estimator proposed by Peng (2001). Indeed, $\hat{\mu}$ has the same asymptotic behavior as the estimator in Peng (2001) since no side information is involved in our empirical likelihood method.

**3. A simulation study and a real application.** For simplicity, our simulation study and real application only deal with the right tail by assuming that the left endpoint of the underlying distribution is finite. Hence we replace $\hat{c}_L(\alpha_L, \mu)$, $\bar{c}_L$ and $X_{n,m}$ by zero and remove (2.7) in our empirical likelihood method described in Section 2.

3.1. *A simulation study.* In this section we conduct a simulation study to investigate the coverage accuracy of our proposed empirical-likelihood-based confidence interval for a mean, and then we compare it with the normal approximation method and the subsample bootstrap method which are given below.

Define

$$\bar{X}_n = \frac{1}{n} \sum_{i=1}^{n} X_i,$$

$$S_n = \left\{ \frac{1}{n} \sum_{i=1}^{n} (X_i - \bar{X}_n)^2 \right\}^{1/2},$$

$$T_n = \sqrt{n}(\bar{X}_n - \mu)/S_n.$$

Without verifying the finiteness of variance, a confidence interval with nominal level $1 - \alpha$ based on normal approximation can be obtained as

$$I_{1-\alpha}^N = (\bar{X}_n - z_{\alpha/2} S_n/\sqrt{n}, \bar{X}_n + z_{\alpha/2} S_n/\sqrt{n}),$$

where $z_{\alpha/2}$ satisfies $P(|N(0,1)| > z_{\alpha/2}) = \alpha$. The subsample bootstrap method was proposed by Hall and LePage (1996) as follows. Conditional on $X_1, \ldots, X_n$, let $X_1^*, \ldots, X_{n_1}^*$ denote independent and identically distributed random variables drawn randomly, with replacement, from $X_1, \ldots, X_n$. Put

$$\bar{X}_{n_1}^* = \frac{1}{n_1} \sum_{i=1}^{n_1} X_i^*,$$

$$S_{n_1}^* = \left\{ \frac{1}{n_1} \sum_{i=1}^{n_1} (X_i^* - \bar{X}_{n_1}^*)^2 \right\}^{1/2},$$

$$T_{n_1}^* = \sqrt{n_1} \frac{(\bar{X}_{n_1}^* - \bar{X}_n)}{S_{n_1}^*}$$

and

$$\hat{x}_{1-\alpha} = \sup\{x : P(|T_{n_1}^*| \leq x | X_1, \ldots, X_n) \leq 1 - \alpha\}.$$

Then a nominal $1 - \alpha$ level confidence interval for $\mu$ is

$$I_{1-\alpha}^* = (\bar{X}_n - \hat{x}_{1-\alpha} S_n/\sqrt{n}, \bar{X}_n + \hat{x}_{1-\alpha} S_n/\sqrt{n}),$$

which has asymptotic coverage probability $1 - \alpha$ under very mild regularity conditions, including $n_1 \to \infty$ and $n_1/n \to 0$ as $n \to \infty$ [see Hall and LePage (1996) for details].

We generated 500 pseudorandom samples of size $n = 1000$ from one of the two distributions: (1) a Burr$(\alpha_R, \beta_R + \alpha_R)$ distribution, given by $F(x) = 1 - (1 + x^{\beta_R})^{-\alpha_R/\beta_R}$, $x > 0$; (2) a Frechet$(\alpha_R)$ distribution, given by $F(x) = \exp\{-x^{-\alpha_R}\}$, $x > 0$. For the subsample bootstrap method we drew 1000 resamples each time.

First we compare our empirical likelihood method with the normal approximation method and the subsample bootstrap method in terms of coverage probability by employing a practical choice of $k = [k^*/\log(k^*)] = n_1$ in view of Theorem 1, where

$$
\begin{aligned}
k^* &= (2^{-1}\alpha_R \beta_R^2 (\beta_R - \alpha_R)^{-3} b_R^{-2} c_R^{2\beta_R/\alpha_R})^{(\alpha_R)/(2\beta_R - \alpha_R)} \\
&\quad \times n^{(2\beta_R - 2\alpha_R)/(2\beta_R - \alpha_R)}
\end{aligned}
\tag{3.1}
$$

minimizes the mean squared error of the Hill estimator [see Hall and Welsh (1985) or de Haan and Peng (1998)]. Here we use the theoretical value of $k^*$ rather than the estimated value, since we investigate the effect of the choice of sample fraction in our next comparison. Coverage probabilities based on these three methods are reported in Table 1. We can conclude from Table 1 that our empirical likelihood method is better than the other two methods when the index $\alpha_R$ is near 2.

Second, we compare our empirical likelihood method with the normal approximation method and the subsample bootstrap method in terms of coverage probability by employing different choices of sample fraction for distributions Frechet(2.0) and Burr(2.0, 4.0); see Tables 2 and 3. These two tables show that the advantage of our empirical likelihood method is robust against the choice of sample fraction.

Third, we investigate the lengths of confidence intervals based on these three methods by employing the practical choice of $k = [k^*/\log(k^*)] = n_1$

Table 1
*Coverage probabilities by employing a practical choice of sample fraction*

|                  | $I_{0.90}$ | $I_{0.90}^*$ | $I_{0.90}^N$ | $I_{0.95}$ | $I_{0.95}^*$ | $I_{0.95}^N$ | $k = n_1$ |
|------------------|-----------|-------------|-------------|-----------|-------------|-------------|-----------|
| Frechet(1.5)     | 0.894     | 0.888       | 0.728       | 0.944     | 0.950       | 0.782       | 37        |
| Frechet(1.8)     | 0.900     | 0.926       | 0.830       | 0.952     | 0.970       | 0.868       | 37        |
| Frechet(2.0)     | 0.910     | 0.936       | 0.850       | 0.954     | 0.980       | 0.902       | 37        |
| Frechet(2.2)     | 0.910     | 0.940       | 0.874       | 0.956     | 0.978       | 0.920       | 37        |
| Frechet(3.0)     | 0.908     | 0.940       | 0.902       | 0.968     | 0.984       | 0.932       | 37        |
| Frechet(5.0)     | 0.914     | 0.922       | 0.896       | 0.968     | 0.980       | 0.946       | 37        |
| Burr(1.5, 3.0)   | 0.886     | 0.906       | 0.728       | 0.946     | 0.962       | 0.784       | 25        |
| Burr(1.8, 3.6)   | 0.894     | 0.938       | 0.834       | 0.954     | 0.980       | 0.868       | 25        |
| Burr(2.0, 4.0)   | 0.896     | 0.946       | 0.854       | 0.952     | 0.982       | 0.900       | 25        |
| Burr(2.2, 4.4)   | 0.910     | 0.950       | 0.874       | 0.956     | 0.984       | 0.924       | 25        |
| Burr(3.0, 6.0)   | 0.906     | 0.944       | 0.900       | 0.960     | 0.988       | 0.936       | 25        |
| Burr(5.0, 10.0)  | 0.910     | 0.936       | 0.904       | 0.966     | 0.976       | 0.942       | 25        |

Note: We report the coverage probabilities for confidence intervals based on our empirical likelihood method, the normal approximation method and the subsample bootstrap method with confidence levels 0.90 and 0.95. We choose $k = [k^*/\log(k^*)] = n_1$, where $k^*$ is defined in (3.1).

TABLE 2
*Coverage probabilities by employing different choices of sample fraction for* Frechet(2.0)

| $k = n_1$ | $I_{0.90}$ | $I^*_{0.90}$ | $I^N_{0.90}$ | $I_{0.95}$ | $I^*_{0.95}$ | $I^N_{0.95}$ |
|---|---|---|---|---|---|---|
| 20 | 0.900 | 0.956 | 0.850 | 0.950 | 0.988 | 0.902 |
| 22 | 0.900 | 0.956 | 0.850 | 0.946 | 0.990 | 0.902 |
| 24 | 0.898 | 0.950 | 0.850 | 0.958 | 0.986 | 0.902 |
| 26 | 0.892 | 0.946 | 0.850 | 0.954 | 0.988 | 0.902 |
| 28 | 0.898 | 0.940 | 0.850 | 0.956 | 0.982 | 0.902 |
| 30 | 0.904 | 0.936 | 0.850 | 0.956 | 0.980 | 0.902 |
| 32 | 0.898 | 0.936 | 0.850 | 0.958 | 0.980 | 0.902 |
| 34 | 0.918 | 0.936 | 0.850 | 0.952 | 0.980 | 0.902 |
| 36 | 0.908 | 0.936 | 0.850 | 0.954 | 0.980 | 0.902 |
| 38 | 0.908 | 0.934 | 0.850 | 0.958 | 0.974 | 0.902 |
| 40 | 0.912 | 0.932 | 0.850 | 0.958 | 0.976 | 0.902 |
| 42 | 0.906 | 0.934 | 0.850 | 0.954 | 0.972 | 0.902 |
| 44 | 0.906 | 0.930 | 0.850 | 0.958 | 0.968 | 0.902 |
| 46 | 0.900 | 0.930 | 0.850 | 0.960 | 0.968 | 0.902 |
| 48 | 0.898 | 0.930 | 0.850 | 0.962 | 0.966 | 0.902 |
| 50 | 0.900 | 0.926 | 0.850 | 0.962 | 0.966 | 0.902 |
| 52 | 0.902 | 0.928 | 0.850 | 0.958 | 0.964 | 0.902 |
| 54 | 0.904 | 0.926 | 0.850 | 0.958 | 0.964 | 0.902 |
| 56 | 0.906 | 0.926 | 0.850 | 0.958 | 0.962 | 0.902 |
| 58 | 0.906 | 0.926 | 0.850 | 0.960 | 0.960 | 0.902 |
| 60 | 0.908 | 0.922 | 0.850 | 0.954 | 0.960 | 0.902 |

Note: We report the coverage probabilities for confidence intervals based on our empirical likelihood method, the normal approximation method and the subsample bootstrap method with confidence levels 0.90 and 0.95. The underlying distribution is taken as Frechet(2.0).

given above. We took confidence level 90% and considered distributions Frechet(1.8), Frechet(2.0), Frechet(2.2), Burr(1.8, 3.6), Burr(2.0, 4.0) and Burr(2.2, 4.4). Since we found that our empirical likelihood method gave a very large value for the right endpoint from time to time, we report in Table 4 the median of left and right endpoints of the confidence intervals based on these three methods. Note that the distributions considered above are skewed to the right and our empirical likelihood method clearly shows this property, that is, with a larger median of endpoints. This also explains why our empirical likelihood method is better than the other two methods in terms of coverage probability when the tail index is near 2. Moreover, our box plots for the endpoints, which are not presented here, show that the symmetric intervals based on the subsample bootstrap method are not good since there are some left endpoints have negative values when $\alpha_R \leq 2$.

Table 3
*Coverage probabilities by employing different choices of sample fraction for* Burr$(2.0, 4.0)$

| $k = n_1$ | $I_{0.90}$ | $I_{0.90}^*$ | $I_{0.90}^N$ | $I_{0.95}$ | $I_{0.95}^*$ | $I_{0.95}^N$ |
|---|---|---|---|---|---|---|
| 20 | 0.904 | 0.948 | 0.854 | 0.952 | 0.986 | 0.900 |
| 22 | 0.902 | 0.948 | 0.854 | 0.946 | 0.984 | 0.900 |
| 24 | 0.896 | 0.948 | 0.854 | 0.954 | 0.986 | 0.900 |
| 26 | 0.894 | 0.946 | 0.854 | 0.952 | 0.982 | 0.900 |
| 28 | 0.900 | 0.938 | 0.854 | 0.954 | 0.978 | 0.900 |
| 30 | 0.906 | 0.938 | 0.854 | 0.956 | 0.976 | 0.900 |
| 32 | 0.896 | 0.936 | 0.854 | 0.956 | 0.978 | 0.900 |
| 34 | 0.916 | 0.936 | 0.854 | 0.952 | 0.976 | 0.900 |
| 36 | 0.908 | 0.932 | 0.854 | 0.956 | 0.978 | 0.900 |
| 38 | 0.910 | 0.932 | 0.854 | 0.958 | 0.974 | 0.900 |
| 40 | 0.916 | 0.932 | 0.854 | 0.958 | 0.970 | 0.900 |
| 42 | 0.908 | 0.934 | 0.854 | 0.952 | 0.970 | 0.900 |
| 44 | 0.908 | 0.932 | 0.854 | 0.960 | 0.970 | 0.900 |
| 46 | 0.900 | 0.930 | 0.854 | 0.958 | 0.968 | 0.900 |
| 48 | 0.898 | 0.926 | 0.854 | 0.962 | 0.964 | 0.900 |
| 50 | 0.904 | 0.930 | 0.854 | 0.960 | 0.966 | 0.900 |
| 52 | 0.900 | 0.930 | 0.854 | 0.958 | 0.964 | 0.900 |
| 54 | 0.900 | 0.926 | 0.854 | 0.958 | 0.962 | 0.900 |
| 56 | 0.906 | 0.922 | 0.854 | 0.960 | 0.962 | 0.900 |
| 58 | 0.906 | 0.922 | 0.854 | 0.958 | 0.960 | 0.900 |
| 60 | 0.912 | 0.922 | 0.854 | 0.956 | 0.962 | 0.900 |

Note: We report the coverage probabilities for confidence intervals based on our empirical likelihood method, the normal approximation method and the subsample bootstrap method with confidence levels 0.90 and 0.95. The underlying distribution is taken as Burr$(2.0, 4.0)$.

3.2. *A real application.* The data set we analyzed consists of 2156 Danish fire losses over one million Danish krone from the years 1980 to 1990 inclusive (see Figure 1). The loss figure is a total loss figure for the event concerned and includes damage to buildings, furnishings and personal property as well as loss of profits. This Danish fire data set was analyzed by McNeil (1997) and Resnick (1996), where the right tail index was confirmed to be between 1 and 2. Here we try to find a confidence interval for the mean value of the Danish fire loss. In Figure 2 we plot the endpoints of the confidence intervals based on the normal approximation method, our empirical likelihood method and the subsample bootstrap method with various choices of sample fraction $k$ and subsample size $n_1$. We observe that the empirical-likelihood-based confidence interval has larger values than the interval constructed from the subsample bootstrap method in most cases. The interval based on normal approximation has a significantly smaller length

<> </>

EMPIRICAL LIKELIHOOD WITH HEAVY TAILS 11

TABLE 4
*Median of endpoints of confidence intervals with confidence level 90%*

|  | MoLE–SBM | MoLE–ELM | MoLE–NAM | MoRE–SBM | MoRE–ELM | MoRE–NAM |
|---|---|---|---|---|---|---|
| Frechet(1.8) | 1.635 | 1.853 | 1.782 | 2.275 | 2.428 | 2.137 |
| Frechet(2.0) | 1.547 | 1.683 | 1.633 | 1.961 | 2.013 | 1.884 |
| Frechet(2.2) | 1.472 | 1.559 | 1.527 | 1.765 | 1.789 | 1.717 |
| Burr(1.8, 3.6) | 1.394 | 1.632 | 1.559 | 2.080 | 2.242 | 1.920 |
| Burr(2.0, 4.0) | 1.333 | 1.471 | 1.428 | 1.770 | 1.816 | 1.686 |
| Burr(2.2, 4.4) | 1.278 | 1.373 | 1.336 | 1.586 | 1.603 | 1.533 |

Note: We report the median of both left and right endpoints of the confidence intervals based on the subsample bootstrap method (SBM), the empirical likelihood method (ELM) and the normal approximation method (NAM) with confidence level 90%. We use MoLE–SBM and MoRE–SBM to denote the median of left and right endpoints of the confidence intervals based on the subsample bootstrap method, respectively; similarly for notation MoLE–ELM, MoRE–ELM, MoLE–NAM and MoRE–NAM.

than the other two types of intervals due to infinite variance. Moreover, our empirical likelihood method clearly demonstrates the skewness of the data to the right and gives asymmetric confidence intervals. As we mentioned in Remark 2, the choice of sample fraction is difficult and very important. We hope to be able to carry on our research on this issue in the future.

**4. Proofs.**

PROOF OF (2.8). Put

$$g_1 = g_1(c_R, c_L, \alpha_R, \alpha_L) = 1 - c_R X_{n,n-k+1}^{-\alpha_R} - c_L(-X_{n,m})^{-\alpha_L},$$

$$g_2 = g_2(\mu, c_R, c_L, \alpha_R, \alpha_L)$$
$$= \mu - c_R \frac{\alpha_R}{\alpha_R - 1} X_{n,n-k+1}^{1-\alpha_R} - c_L \frac{\alpha_L}{\alpha_L - 1}(-X_{n,m})^{1-\alpha_L},$$

$$H(\mu, c_R, c_L, \alpha_R, \alpha_L, \lambda_1, \lambda_2, p_i)$$
$$= k \log \alpha_R + k \log c_R - (\alpha_R + 1) \sum_{\bar{\delta}_i = 1} \log X_i$$
$$+ m \log \alpha_L + m \log c_L - (\alpha_L + 1) \sum_{\delta_i = 1} \log X_i$$
$$+ \sum_{\bar{\delta}_i = \delta_i = 0} \log p_i - \lambda_1 \left\{ \sum_{\bar{\delta}_i = \delta_i = 0} p_i - g_1 \right\} - \lambda_2 \sum_{\bar{\delta}_i = \delta_i = 0} p_i w_i.$$

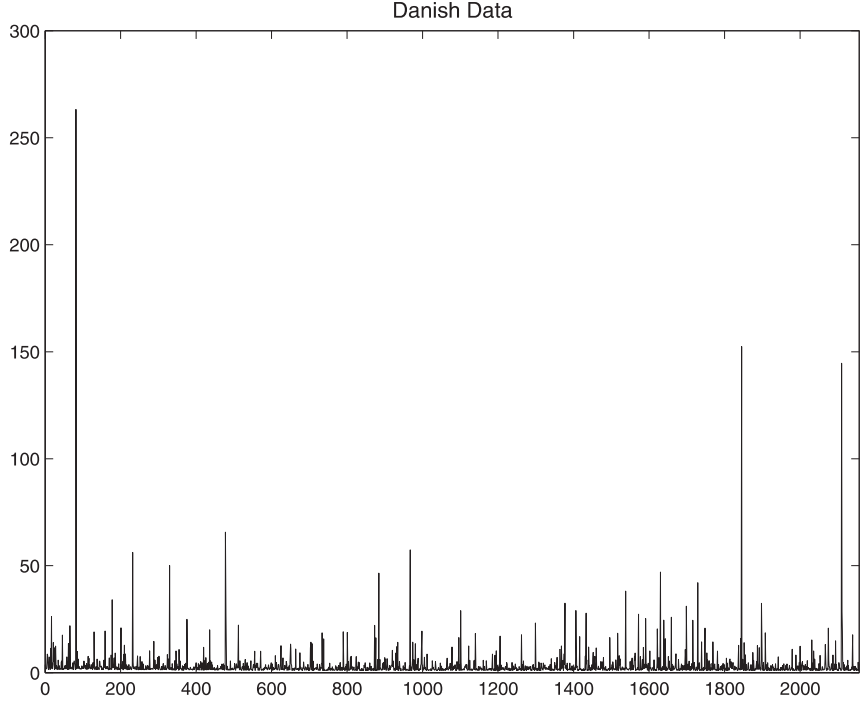

Fig. 1. *Danish fire loss data, which consists of* 2156 *losses over one million Danish krone for the years* 1980 *to* 1990 *inclusive.*

By the method of Lagrange multipliers we have

$$
\begin{aligned}
(4.1) \quad & \frac{k}{c_R} - \lambda_1 X_{n,n-k+1}^{-\alpha_R} \\
& - \lambda_2 \left\{ \frac{\alpha_R}{\alpha_R - 1} X_{n,n-k+1}^{1-\alpha_R} - g_1^{-1} g_2 X_{n,n-k+1}^{-\alpha_R} \right\} = 0,
\end{aligned}
$$

$$
\begin{aligned}
(4.2) \quad & \frac{m}{c_L} - \lambda_1 (-X_{n,m})^{-\alpha_L} \\
& - \lambda_2 \left\{ \frac{\alpha_L}{\alpha_L - 1} (-X_{n,m})^{1-\alpha_L} - g_1^{-1} g_2 (-X_{n,m})^{-\alpha_L} \right\} = 0,
\end{aligned}
$$

$$
\begin{aligned}
(4.3) \quad & \frac{k}{\alpha_R} - \sum_{\bar{\delta}_i = 1} \log X_i + \lambda_1 c_R X_{n,n-k+1}^{-\alpha_R} \log X_{n,n-k+1} \\
& + \lambda_2 \Bigg\{ \frac{c_R}{(\alpha_R - 1)^2} X_{n,n-k+1}^{1-\alpha_R} + \frac{c_R \alpha_R}{\alpha_R - 1} X_{n,n-k+1}^{1-\alpha_R} \log X_{n,n-k+1} \\
& \qquad - g_1^{-1} g_2 c_R X_{n,n-k+1}^{-\alpha_R} \log X_{n,n-k+1} \Bigg\} = 0,
\end{aligned}
$$

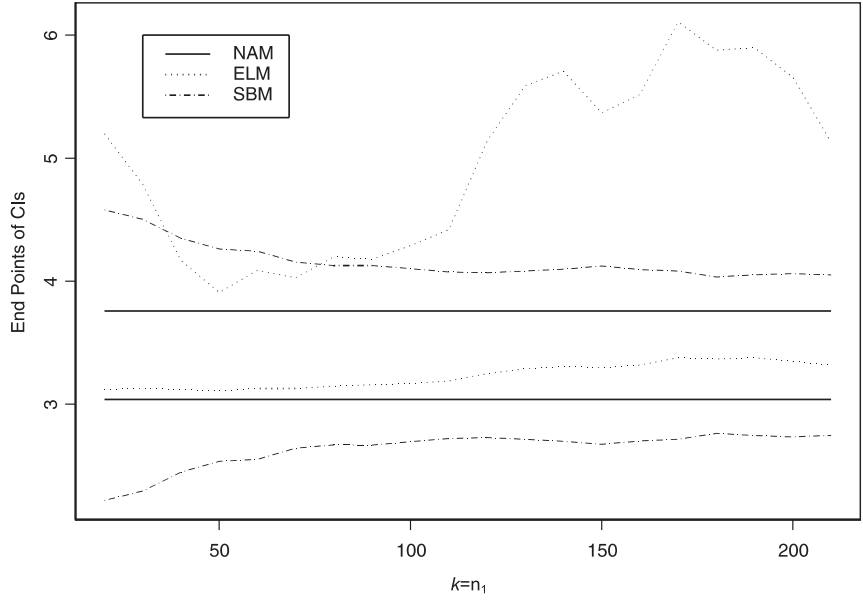

FIG. 2. *Confidence intervals for the mean value of Danish fire loss data. The endpoints of confidence intervals with level 95% based on our empirical likelihood method, the normal approximation method and the subsample bootstrap method are plotted against $k = n_1 = 20, 30, 40, \ldots, 210$. Here NAM, ELM and SBM denote the normal approximation method, our empirical likelihood method and the subsample bootstrap method, respectively.*

$$\begin{aligned}
& \frac{m}{\alpha_L} - \sum_{\delta_i=1} \log X_i + \lambda_1 c_L (-X_{n,m})^{-\alpha_L} \log(-X_{n,m}) \\
& \quad + \lambda_2 \bigg\{ \frac{c_L}{(\alpha_L - 1)^2} (-X_{n,m})^{1-\alpha_L} \\
& \quad\quad + \frac{c_L \alpha_L}{\alpha_L - 1} (-X_{n,m})^{1-\alpha_L} \log(-X_{n,m}) \\
& \quad\quad - g_1^{-1} g_2 c_L (-X_{n,m})^{-\alpha_L} \log(-X_{n,m}) \bigg\} = 0,
\end{aligned}$$
(4.4)

$$p_i = \{\lambda_1 + \lambda_2 w_i\}^{-1}, \qquad \text{if } \bar{\delta}_i = \delta_i = 0, \tag{4.5}$$

$$\sum_{\bar{\delta}_i = \delta_i = 0} \frac{1}{\lambda_1 + \lambda_2 w_i} = g_1, \tag{4.6}$$

$$\sum_{\bar{\delta}_i = \delta_i = 0} \frac{w_i}{\lambda_1 + \lambda_2 w_i} = 0. \tag{4.7}$$

By (4.1) and (4.3),

$$\text{(4.8)} \quad \frac{k}{\alpha_R} - \sum_{\bar{\delta}_i=1} \log \frac{X_i}{X_{n,n-k+1}} + \frac{\lambda_2 c_R}{(\alpha_R-1)^2} X_{n,n-k+1}^{1-\alpha_R} = 0.$$

By (4.2) and (4.4),

$$\text{(4.9)} \quad \frac{m}{\alpha_L} - \sum_{\delta_i=1} \log \frac{X_i}{X_{n,m}} + \frac{\lambda_2 c_L}{(\alpha_L-1)^2} (-X_{n,m})^{1-\alpha_L} = 0.$$

By (4.5)–(4.7),

$$\text{(4.10)} \quad \begin{aligned} \lambda_1 &= \frac{n-m-k}{g_1} \\ &= (n-m-k)\{1 - c_R X_{n,n-k+1}^{-\alpha_R} - c_L(-X_{n,m})^{-\alpha_L}\}^{-1}. \end{aligned}$$

By (4.1), (4.2) and (4.10),

$$\text{(4.11)} \quad \lambda_2 g_1^{-1} g_2 = \lambda_2 \mu + \lambda_1 - n.$$

By (4.1) and (4.11),

$$\text{(4.12)} \quad \frac{k}{c_R} - \frac{\lambda_2 \alpha_R}{\alpha_R - 1} X_{n,n-k+1}^{1-\alpha_R} + \lambda_2 \mu X_{n,n-k+1}^{-\alpha_R} - n X_{n,n-k+1}^{-\alpha_R} = 0.$$

By (4.2) and (4.11),

$$\text{(4.13)} \quad \frac{m}{c_L} - \frac{\lambda_2 \alpha_L}{\alpha_L - 1} (-X_{n,m})^{1-\alpha_L} + \lambda_2 \mu (-X_{n,m})^{-\alpha_L} - n(-X_{n,m})^{-\alpha_L} = 0.$$

By (4.8) and (4.12),

$$\text{(4.14)} \quad \begin{aligned} \lambda_2 = \Bigg\{ & kn^{-1}(\alpha_R-1)^{-2} X_{n,n-k+1} \bigg[\sum_{\bar{\delta}_i=1} \log \frac{X_i}{X_{n,n-k+1}} - \frac{k}{\alpha_R}\bigg]^{-1} \\ & - n^{-1}\alpha_R(\alpha_R-1)^{-1} X_{n,n-k+1} + n^{-1}\mu \Bigg\}^{-1} \end{aligned}$$

and

$$\text{(4.15)} \quad \begin{aligned} c_R = & \frac{k}{n} X_{n,n-k+1}^{\alpha_R} \\ & - \bigg\{\sum_{\bar{\delta}_i=1} \log \frac{X_i}{X_{n,n-k+1}} - \frac{k}{\alpha_R}\bigg\} \\ & \times \{\alpha_R(\alpha_R-1)n^{-1} X_{n,n-k+1}^{\alpha_R} - \mu(\alpha_R-1)^2 n^{-1} X_{n,n-k+1}^{\alpha_R-1}\}. \end{aligned}$$

By (4.9) and (4.13),

$$
\begin{aligned}
\lambda_2 = \bigg\{ & mn^{-1}(\alpha_L-1)^{-2}(-X_{n,m})\bigg[\sum_{\delta_i=1}\log\frac{X_i}{X_{n,m}}-\frac{m}{\alpha_L}\bigg]^{-1} \\
& - n^{-1}\alpha_L(\alpha_L-1)^{-1}(-X_{n,m})+n^{-1}\mu\bigg\}^{-1}
\end{aligned}
\tag{4.16}
$$

and

$$
\begin{aligned}
c_L = &\; \frac{m}{n}(-X_{n,m})^{\alpha_L} \\
& -\bigg\{\sum_{\delta_i=1}\log\frac{X_i}{X_{n,m}}-\frac{m}{\alpha_L}\bigg\} \\
& \times\{\alpha_L(\alpha_L-1)n^{-1}(-X_{n,m})^{\alpha_L}-\mu(\alpha_L-1)^2 n^{-1}(-X_{n,m})^{\alpha_L-1}\}.
\end{aligned}
\tag{4.17}
$$

So (2.6) follows from (4.7), (4.10), (4.14), (4.15) and (4.17). Similarly (2.7) follows from (4.7), (4.10) and (4.15)–(4.17). Hence (2.8). $\square$

Before we prove Theorem 1, we need some notation and several lemmas. Let $F^-$ denote the inverse function of $F$. Then (2.1) implies that

$$
\begin{aligned}
F^-(1-x) &= c_R^{1/\alpha_R} x^{-1/\alpha_R}\{1+d_R x^{\rho_R}+o(x^{\rho_R})\}, \\
F^-(x) &= -c_L^{1/\alpha_L} x^{-1/\alpha_L}\{1+d_L x^{\rho_L}+o(x^{\rho_L})\}
\end{aligned}
\tag{4.18}
$$

as $x \to 0$, where $\rho_R = \beta_R/\alpha_R$, $d_R = b_R c_R^{-\rho_R}/\alpha_R$, $\rho_L = \beta_L/\alpha_L$ and $d_L = b_L c_L^{-\rho_L}/\alpha_L$.

Let $U_1, U_2, \ldots$ be i.i.d. random variables with a uniform distribution on $(0,1)$ and let $U_{n,1} \leq \cdots \leq U_{n,n}$ denote the order statistics of $U_1, \ldots, U_n$. Define $G_n(u) = \frac{1}{n}\sum_{i=1}^n I(U_i \leq u)$. Then from Csörgő, Csörgő, Horváth and Mason (1986) there exists a sequence of Brownian bridges $B_n(u)$, $0 \leq u \leq 1$, $n = 1, 2, \ldots$, such that

$$
\sup_{1/n \leq u \leq 1-1/n} n^\mu \frac{|\sqrt{n}(G_n(u)-u)-B_n(u)|}{(u(1-u))^{1/2-\mu}} = O_p(1)
\tag{4.19}
$$

as $n \to \infty$, where $\mu$ is any fixed number such that $\mu \in [0, 1/4)$. Without loss of generality we assume $X_{n,j} = F^-(U_{n,j})$, $j = 1, 2, \ldots, n$. We use $\mu_0$, $c_{R0}$, $c_{L0}$, $\alpha_{R0}$ and $\alpha_{L0}$ to denote the true parameters under model (1.1).

LEMMA 1. *Put*

$$
\sigma^2(s,t) = \int_s^{1-t}\int_s^{1-t}(u\wedge v - uv)\,dF^-(u)\,dF^-(v).
$$

*Suppose* (2.1) *holds. Then*

$$\sigma^2\left(\frac{m}{n},\frac{k}{n}\right) \sim \begin{cases} \int_{m/n}^{1-k/n}[F^-(s)]^2\,ds, & \text{if } \alpha_R \geq 2, \alpha_L \geq 2, \\ \left(\frac{m}{n}\right)^{1-2/\alpha_L}\left(c_L^2 + \frac{1}{2-\alpha_L}\right), & \text{if } \alpha_R > 2, \alpha_L < 2, \\ \left(\frac{m}{n}\right)^{1-2/\alpha_L}\left(c_L^2 + \frac{1}{2-\alpha_L}\right) + \int_{1/2}^{1-k/n}[F^-(s)]^2\,ds, & \\ & \text{if } \alpha_R = 2, \alpha_L < 2, \\ \left(\frac{k}{n}\right)^{1-2/\alpha_R}\left(c_R^2 + \frac{1}{2-\alpha_R}\right), & \text{if } \alpha_R < 2, \alpha_L > 2, \\ \left(\frac{k}{n}\right)^{1-2/\alpha_R}\left(c_R^2 + \frac{1}{2-\alpha_R}\right) + \int_{m/n}^{1/2}[F^-(s)]^2\,ds, & \\ & \text{if } \alpha_R < 2, \alpha_L = 2, \\ \left(\frac{k}{n}\right)^{1-2/\alpha_R}\left(c_R^2 + \frac{1}{2-\alpha_R}\right) + \left(\frac{m}{n}\right)^{1-2/\alpha_L}\left(c_L^2 + \frac{1}{2-\alpha_L}\right), & \\ & \text{if } \alpha_R < 2, \alpha_L < 2. \end{cases}$$

*Moreover, if* $B_0 n^{l_0} \leq m, k \leq B_1 n^{l_1}$ *for some* $B_0 > 0$, $B_1 > 0$, $0 < l_0 \leq l_1 < 1$, *then*

$$\int_{1/2}^{1-k/n}[F^-(s)]^2\,ds/(m/n)^{1-2/\alpha_L} \to 0 \quad \text{in case} \quad \alpha_R = 2, \alpha_L < 2,$$

$$\int_{m/n}^{1/2}[F^-(s)]^2\,ds/(k/n)^{1-2/\alpha_R} \to 0 \quad \text{in case} \quad \alpha_R < 2, \alpha_L = 2.$$

PROOF.  It is easy to check that

$$\sigma^2\left(\frac{m}{n},\frac{k}{n}\right) = \frac{m}{n}\left[F^-\left(\frac{m}{n}\right)\right]^2 - \left(\frac{m}{n}\right)^2\left[F^-\left(\frac{m}{n}\right)\right]^2$$

$$-2\frac{m}{n}F^-\left(\frac{m}{n}\right)\int_{m/n}^{1-k/n} F^-(s)\,ds$$

$$+\frac{k}{n}\left[F^-\left(1-\frac{k}{n}\right)\right]^2 - \left(\frac{k}{n}\right)^2\left[F^-\left(1-\frac{k}{n}\right)\right]^2$$

$$-2\frac{k}{n}F^-\left(1-\frac{k}{n}\right)\int_{m/n}^{1-k/n} F^-(s)\,ds$$

$$+\int_{m/n}^{1-k/n}\left[F^-(s)\right]^2 ds - \left[\int_{m/m}^{1-k/n} F^-(s)\,ds\right]^2.$$

The rest of the proof follows by applying Potter's bounds [see Geluk and de Haan (1987)] to (4.18) in a similar way to the proof of Lemma 2.5 in Csörgő, Haeusler and Mason (1988).  □

LEMMA 2. *Suppose the conditions in Theorem 1 hold. Then the following approximations are true:*

$$\sqrt{k}\{\bar{\alpha}_R - \alpha_{R0}\}$$
$$= -\alpha_{R0}\sqrt{\frac{n}{k}}B_n\left(1 - \frac{k}{n}\right) + \alpha_{R0}\sqrt{\frac{n}{k}}\int_{1-k/n}^{1}\frac{B_n(s)}{1-s}ds + o_p(1),$$

$$\sqrt{m}\{\bar{\alpha}_L - \alpha_{L0}\}$$
$$= \alpha_{L0}\sqrt{\frac{n}{m}}B_n\left(\frac{m}{n}\right) - \alpha_{L0}\sqrt{\frac{n}{m}}\int_{0}^{m/n}\frac{B_n(s)}{s}ds + o_p(1),$$

$$\sqrt{k}\left\{\frac{\bar{c}_R}{c_{R0}} - 1\right\} = \sqrt{\frac{n}{k}}B_n\left(1 - \frac{k}{n}\right) + o_p(1),$$

$$\sqrt{k}\left\{\frac{\bar{c}_L}{c_{L0}} - 1\right\} = -\sqrt{\frac{n}{m}}B_n\left(\frac{m}{n}\right) + o_p(1),$$

$$\sqrt{k}\left\{\frac{\hat{c}_R(\alpha_{R0},\mu_0)}{c_{R0}} - 1\right\}$$
$$= \sqrt{\frac{n}{k}}B_n\left(1 - \frac{k}{n}\right)$$
$$- \alpha_{R0}^{-1}\left\{\sqrt{\frac{n}{k}}B_n\left(1 - \frac{k}{n}\right) - \sqrt{\frac{n}{k}}\int_{1-k/n}^{1}\frac{B_n(s)}{1-s}ds\right\}$$
$$\times\left\{\alpha_{R0}(\alpha_{R0} - 1) - \mu_0(\alpha_{R0} - 1)^2 c_{R0}^{-1/\alpha_{R0}}\left(\frac{k}{n}\right)^{1/\alpha_{R0}}\right\} + o_p(1),$$

$$\sqrt{m}\left\{\frac{\hat{c}_L(\alpha_{L0},\mu_0)}{c_{L0}} - 1\right\}$$
$$= -\sqrt{\frac{n}{m}}B_n\left(\frac{m}{n}\right)$$
$$- \alpha_{L0}^{-1}\left\{-\sqrt{\frac{n}{m}}B_n\left(\frac{m}{n}\right) + \sqrt{\frac{n}{m}}\int_{0}^{m/n}\frac{B_n(s)}{s}ds\right\}$$
$$\times\left\{\alpha_{L0}(\alpha_{L0} - 1) - \mu_0(\alpha_{L0} - 1)^2 c_{L0}^{-1/\alpha_{L0}}\left(\frac{m}{n}\right)^{1/\alpha_{L0}}\right\} + o_p(1),$$

$$\hat{\lambda}_2\frac{(\alpha_{R0},\mu_0)}{n} = o_p\left(\left(\frac{k}{n}\right)^{1/\alpha_{R0}}k^{-1/2}\right),$$

$$\bar{\lambda}_2\frac{(\alpha_{L0},\mu_0)}{n} = o_p\left(\left(\frac{m}{n}\right)^{1/\alpha_{L0}}m^{-1/2}\right),$$

$$\hat{\lambda}_1 \frac{(\alpha_{R0}, \alpha_{L0}, \mu_0)}{n} = 1 + o_p(1),$$

$$\frac{\sqrt{n}}{\sigma(m/n, k/n)} \left\{ \frac{1}{n} \sum_{\bar{\delta}_i = \delta_i = 0} X_i - \int_{m/n}^{1-k/n} F^-(s)\, ds \right\}$$

$$= -\int_{m/n}^{1-k/n} \frac{B_n(s)\, dF^-(s)}{\sigma(m/k, k/n)} + o_p(1)$$

$$= \frac{\sqrt{n}}{\sigma(m/n, k/n)} \left\{ \frac{1}{n} \sum_{\bar{\delta}_i = \delta_i = 0} w_i(\mu_0, c_{R0}, c_{L0}, \alpha_{R0}, \alpha_{L0}) \right\} + o_p(1).$$

PROOF. Using Theorems 2.3 and 2.4 of Csörgő and Mason (1985), we have

(4.20)
$$\sqrt{k} \left\{ -\frac{1}{k} \sum_{i=1}^{k} \log \frac{1 - U_{n,n-i+1}}{1 - U_{n,n-k+1}} - 1 \right\}$$
$$= \sqrt{\frac{n}{k}} B_n\left(1 - \frac{k}{n}\right) - \sqrt{\frac{n}{k}} \int_{1-k/n}^{1} \frac{B_n(s)}{1-s}\, ds + o_p(1),$$

(4.21)
$$\sqrt{k} \left\{ -\log(1 - U_{n,n-k+1}) + \log \frac{k+1}{n} \right\}$$
$$= \sqrt{\frac{n}{k}} B_n\left(1 - \frac{k}{n}\right) + o_p(1).$$

Similarly

(4.22)
$$\sqrt{m} \left\{ -\frac{1}{m} \sum_{i=1}^{m} \log \frac{U_{n,i}}{U_{n,m}} - 1 \right\}$$
$$= -\sqrt{\frac{n}{m}} B_n\left(\frac{m}{n}\right) + \sqrt{\frac{n}{m}} \int_{0}^{m/n} \frac{B_n(s)}{s}\, ds + o_p(1),$$

(4.23)
$$\sqrt{m} \left\{ -\log U_{n,m} + \log \frac{m}{n} \right\}$$
$$= -\sqrt{\frac{n}{m}} B_n\left(\frac{m}{n}\right) + o_p(1).$$

Hence the lemma follows from (4.18)–(4.23) and Lemma 1. □

LEMMA 3. *Suppose the conditions in Theorem 1 hold. Then there exist solutions $(\hat{\alpha}_R(\mu_0), \hat{\alpha}_L(\mu_0))$ to (2.6) and (2.7) such that $\hat{\alpha}_R(\mu_0) - \alpha_{R0} = O_p(\delta_q)$ and $\hat{\alpha}_L(\mu_0) - \alpha_{L0} = O_p(\delta_q)$, where $\delta_q = O(k^{-q} + m^{-q})$ for any $q \in [1/3, 1/2)$.*

PROOF. By Taylor expansion we can easily show that

$$\hat{c}_R(\alpha_{R0} \pm \delta_q, \mu_0) = c_{R0} + o_p(\delta_q),$$
$$\hat{c}_L(\alpha_{L0} \pm \delta_q, \mu_0) = c_{L0} + o_p(\delta_q),$$
$$\hat{\lambda}_2(\alpha_{R0} \pm \delta_q, \mu_0) = n(\alpha_{R0} - 1)^2 \alpha_{R0}^{-2} X_{n,n-k+1}^{-1} \{\pm \delta_q + O_p(k^{-1/2})\},$$
$$\bar{\lambda}_2(\alpha_{L0} \pm \delta_q, \mu_0) = n(\alpha_{L0} - 1)^2 \alpha_{L0}^{-2} (-X_{n,m})^{-1} \{\pm \delta_q + O_p(m^{-1/2})\},$$

$$\hat{h}(\alpha_{R0} \pm \delta_q, \alpha_L, \mu_0)$$
$$= \left\{ -(\alpha_{R0} - 1)^{-2} X_{n,n-k+1}^{1-\alpha_{R0}} \right.$$
$$\left. - \left[ \frac{1}{n} \sum_{\bar{\delta}_i = \delta_i = 0} w_i(\mu_0, c_{R0}, c_L, \alpha_{R0}, \alpha_L) \right] (\alpha_{R0} - 1)^2 \alpha_{R0}^{-1} X_{n,n-k+1}^{-1} \right\}$$
$$\times \{\pm \delta_q + O_p(k^{-1/2})\},$$

$$\hat{h}(\alpha_{R0}, \alpha_{L0} \pm \delta_q, \mu_0) = -(\alpha_{L0} - 1)^{-2}(-X_{n,m})^{1-\alpha_{L0}}\{\pm \delta_q + O_p(m^{-1/2})\},$$
$$\bar{h}(\alpha_{R_0} \pm \delta_q, \alpha_L, \mu_0) = -(\alpha_{R0} - 1)^{-2} X_{n,n-k+1}^{1-\alpha_{R0}}\{\pm \delta_q + O_p(k^{-1/2})\},$$

$$\bar{h}(\alpha_R, \alpha_{L0} \pm \delta_q, \mu_0)$$
$$= \left\{ -(\alpha_{L0} - 1)^{-2}(-X_{n,m})^{1-\alpha_{L0}} \right.$$
$$\left. - \left[ \frac{1}{n} \sum_{\bar{\delta}_i = \delta_i = 0} w_i(\mu_0, c_R, c_{L0}, \alpha_R, \alpha_{L0}) \right] (\alpha_{L0} - 1)^2 \alpha_{L0}^{-1} (-X_{n,m})^{-1} \right\}$$
$$\times \{\pm \delta_q + O_p(m^{-1/2})\}.$$

Hence the lemma follows from the above expansions and

$$\frac{X_{n,m}}{F^-(m/n)} - 1 = O_p\left(\frac{1}{\sqrt{m}}\right),$$
$$\frac{X_{n,n-k+1}}{F^-(1-k/n)} - 1 = O_p\left(\frac{1}{\sqrt{k}}\right). \qquad \square$$

PROOF OF THEOREM 1. Let $\eta = (c_R, c_L, \alpha_R, \alpha_L, p_i)$. The key idea in the proof is to expand $l_0(\hat{\eta})$ around $\bar{\eta}$ and to derive the convergence rates for $\hat{\alpha}_R - \bar{\alpha}_R$, $\hat{\alpha}_L - \bar{\alpha}_L$, $\hat{c}_R - \bar{c}_R$ and $\hat{c}_L - \bar{c}_L$.

Using Lemma 3 we can easily prove that

$$\hat{\lambda}_2(\hat{\alpha}_R, \mu_0) = nX_{n,n-k+1}^{-1}\left\{(\alpha_{R0}-1)^2\left(\frac{1}{k}\sum_{\bar{\delta}_i=1}\log\frac{X_i}{X_{n,n-k+1}} - \alpha_{R0}^{-1}\right)\right.$$

$$\left. + (\alpha_{R0}-1)^2\alpha_{R0}^{-2}(\hat{\alpha}_R - \alpha_{R0}) + o_p\left(\frac{1}{\sqrt{k}}\right)\right\},$$

$$\bar{\lambda}_2(\hat{\alpha}_L, \mu_0) = n(-X_{n,m})^{-1}\left\{(\alpha_{L0}-1)^2\left(\frac{1}{m}\sum_{\delta_i=1}\log\frac{X_i}{X_{n,m}} - \alpha_{L0}^{-1}\right)\right.$$

$$\left. + (\alpha_{L0}-1)^2\alpha_{L0}^{-2}(\hat{\alpha}_L - \alpha_{L0}) + o_p\left(\frac{1}{\sqrt{m}}\right)\right\}.$$

It follows from (2.6) and (2.7) that

$$\hat{\lambda}_2(\hat{\alpha}_R, \mu_0) = \frac{\hat{\lambda}_1(\hat{\alpha}_R, \hat{\alpha}_L, \mu_0)\sum_{\bar{\delta}_i=\delta_i=0}w_i(\mu_0, \hat{c}_R, \hat{c}_L, \hat{\alpha}_R, \hat{\alpha}_L)}{\sum_{\bar{\delta}_i=\delta_i=0}w_i^2(\mu_0, \hat{c}_R, \hat{c}_L, \hat{\alpha}_R, \hat{\alpha}_L)}$$

$$\times \{1 + o_p(1)\}$$

$$= n\sigma^{-2}\left(\frac{m}{n}, \frac{k}{n}\right)\left\{\frac{1}{n}\sum_{\bar{\delta}_i=\delta_i=0}w_i(\mu_0, c_{R0}, c_{L0}, \alpha_{R0}, \alpha_{L0})\right\}$$

(4.24)

$$\times \{1 + o_p(1)\},$$

$$\bar{\lambda}_2(\hat{\alpha}_L, \mu_0) = n\sigma^{-2}\left(\frac{m}{n}, \frac{k}{n}\right)\left\{\frac{1}{n}\sum_{\bar{\delta}_i=\delta_i=0}w_i(\mu_0, c_{R0}, c_{L0}, \alpha_{R0}, \alpha_{L0})\right\}$$

$$\times \{1 + o_p(1)\}.$$

Hence

$$\hat{\alpha}_R - \alpha_{R0} = -\alpha_{R0}^2\left\{\frac{1}{k}\sum_{\bar{\delta}_i=1}\log\frac{X_i}{X_{n,n-k+1}} - \alpha_{R0}^{-1}\right\}$$

$$+ X_{n,n-k+1}(\alpha_{R0}-1)^{-2}\alpha_{R0}^2\sigma^{-2}\left(\frac{m}{n}, \frac{k}{n}\right)$$

(4.25)

$$\times \frac{1}{n}\sum_{\bar{\delta}_i=\delta_i=0}w_i(\mu_0, c_{R0}, c_{L0}, \alpha_{R0}, \alpha_{L0})$$

$$+ o_p\left(\frac{1}{\sqrt{k}}\right)$$

and

$$\hat{\alpha}_L - \alpha_{L0} = -\alpha_{L0}^2 \bigg\{ \frac{1}{m} \sum_{\bar{\delta}_i=1} \log \frac{X_i}{X_{n,m}} - \alpha_{L0}^{-1} \bigg\}$$

$$- X_{n,m}(\alpha_{L0} - 1)^{-2} \alpha_{L0}^2 \sigma^{-2}\left(\frac{m}{n}, \frac{k}{n}\right)$$

(4.26)

$$\times \frac{1}{n} \sum_{\bar{\delta}_i = \delta_i = 0} w_i(\mu_0, c_{R0}, c_{L0}, \alpha_{R0}, \alpha_{L0})$$

$$+ o_p\left(\frac{1}{\sqrt{m}}\right).$$

Furthermore

(4.27) $\quad \dfrac{\hat{c}_R(\hat{\alpha}_R, \mu_0)}{c_{R0}} = \dfrac{\hat{c}_R(\alpha_{R0}, \mu_0)}{c_{R0}} - \alpha_{R0}^{-1}(\alpha_{R0} - 1)(\hat{\alpha}_R - \alpha_{R0})(1 + o_p(1))$

and

(4.28) $\quad \dfrac{\hat{c}_L(\hat{\alpha}_L, \mu_0)}{c_{L0}} = \dfrac{\hat{c}_L(\alpha_{L0}, \mu_0)}{c_{L0}} - \alpha_{L0}^{-1}(\alpha_{L0} - 1)(\hat{\alpha}_L - \alpha_{L0})(1 + o_p(1)).$

Expanding $l_0(\hat{\eta})$ at $\bar{\eta}$ by noting that $\bar{\eta}$ is the maximum likelihood estimate, we have

$$l(\mu_0) = -2 \bigg\{ \frac{1}{2} \frac{\partial^2 l_0}{\partial \alpha_R^2}(\bar{c}_R, \bar{c}_L, \bar{\alpha}_R, \bar{\alpha}_L, \bar{p}_i)[\hat{\alpha}_R - \bar{\alpha}_R]^2$$

$$+ \frac{1}{2} \frac{\partial^2 l_0}{\partial \alpha_L^2}(\bar{c}_R, \bar{c}_L, \bar{\alpha}_R, \bar{\alpha}_L, \bar{p}_i)[\hat{\alpha}_L - \bar{\alpha}_L]^2$$

$$+ \frac{1}{2} \frac{\partial^2 l_0}{\partial c_R^2}(\bar{c}_R, \bar{c}_L, \bar{\alpha}_R, \bar{\alpha}_L, \bar{p}_i)[\hat{c}_R - \bar{c}_R]^2$$

$$+ \frac{1}{2} \frac{\partial^2 l_0}{\partial c_L^2}(\bar{c}_R, \bar{c}_L, \bar{\alpha}_R, \bar{\alpha}_L, \bar{p}_i)[\hat{c}_L - \bar{c}_L]^2$$

$$- \frac{n}{2\sigma^2(m/n, k/n)} \bigg[ \frac{1}{n} \sum_{\bar{\delta}_i = \delta_i = 0} w_i(\mu_0, \hat{c}_R, \hat{c}_L, \hat{\alpha}_R, \hat{\alpha}_L) \bigg]^2 + o_p(1) \bigg\}$$

$$= -2 \bigg\{ -2^{-1} k \alpha_{R0}^{-2}[\hat{\alpha}_R - \bar{\alpha}_R]^2 - 2^{-1} k \alpha_{L0}^{-2}[\hat{\alpha}_L - \bar{\alpha}_L]^2$$

$$- 2^{-1} k c_{R0}^{-2}[\hat{c}_R - \bar{c}_R]^2 - 2^{-1} k c_{L0}^{-2}[\hat{c}_L - \bar{c}_L]^2$$

$$- \frac{n}{2\sigma^2(m/n, k/n)} \bigg[ \frac{1}{n} \sum_{\bar{\delta}_i = \delta_i = 0} w_i(\mu_0, c_{R0}, c_{L0}, \alpha_{R0}, \alpha_{L0}) \bigg]^2 + o_p(1) \bigg\}.$$

It follows from Lemma 2 and (4.24)–(4.28) that

$$k[\hat{\alpha}_R - \bar{\alpha}_R]^2 = o_p(1),$$
$$k[\hat{\alpha}_L - \bar{\alpha}_L]^2 = o_p(1),$$
$$k[\hat{c}_R - \bar{c}_R]^2 = o_p(1),$$
$$k[\hat{c}_L - \bar{c}_L]^2 = o_p(1)$$

and

$$\frac{\sqrt{n}}{\sigma(m/n, k/n)} \left\{ \frac{1}{n} \sum_{\bar{\delta}_i = \delta_i = 0} w_i(\mu_0, c_{R0}, c_{L0}, \alpha_{R0}, \alpha_{L0}) \right\} \xrightarrow{d} N(0, 1).$$

Hence the asymptotic limit of $l(\mu_0)$ is $\chi^2(1)$. $\square$

**Acknowledgments.** I thank Professor Art Owen for discussions and Dr. A. J. McNeil for providing the Danish fire loss data. Support via the Summer New Researcher Program at Stanford University is greatly appreciated. I also thank an Associate Editor and two referees for helpful comments.

## REFERENCES

Chen, S. and Hall, P. (1993). Smoothed empirical likelihood confidence intervals for quantiles. *Ann. Statist.* **21** 1166–1181. MR1241263

Csörgő, M., Csörgő, S., Horváth, L. and Mason, D. M. (1986). Weighted empirical and quantile processes. *Ann. Probab.* **14** 31–85. MR815960

Csörgő, S., Haeusler, E. and Mason, D. M. (1988). The asymptotic distribution of trimmed sums. *Ann. Probab.* **16** 672–699. MR929070

Csörgő, S. and Mason, D. M. (1985). Central limit theorems for sums of extreme values. *Math. Proc. Cambridge Philos. Soc.* **98** 547–558. MR803614

de Haan, L. and Peng, L. (1998). Comparison of tail index estimators. *Statist. Neerlandica* **52** 60–70. MR1615558

de Haan, L. and Stadtmüller, U. (1996). Generalized regular variation of second order. *J. Austral. Math. Soc. Ser. A* **61** 381–395. MR1420345

DiCiccio, T. J., Hall, P. and Romano, J. P. (1991). Empirical likelihood is Bartlett-correctable. *Ann. Statist.* **19** 1053–1061. MR1105861

Embrechts, P., Klüppelberg, C. and Mikosch, T. (1997). *Modelling Extremal Events for Insurance and Finance.* Springer, New York. MR1458613

Feuerverger, A. and Hall, P. (1999). Estimating a tail exponent by modelling departure from a Pareto distribution. *Ann. Statist.* **27** 760–781. MR1714709

Geluk, J. and de Haan, L. (1987). *Regular Variation, Extensions and Tauberian Theorems.* CWI Tract 40. Centrum voor Wiskunde en Informatica, Amsterdam. MR906871

Hall, P. and Jing, B. (1998). Comparison of bootstrap and asymptotic approximations to the distribution of a heavy-tailed mean. *Statist. Sinica* **8** 887–906. MR1651514

Hall, P. and La Scala, B. (1990). Methodology and algorithms of empirical likelihood. *Internat. Statist. Rev.* **58** 109–127.

Hall, P. and LePage, R. (1996). On bootstrap estimation of the distribution of the Studentized mean. *Ann. Inst. Statist. Math.* **48** 403–421. MR1424772

Hall, P. and Welsh, A. H. (1985). Adaptive estimates of parameters of regular variation. *Ann. Statist.* **13** 331–341. MR773171

Hill, B. M. (1975). A simple general approach to inference about the tail of a distribution. *Ann. Statist.* **3** 1163–1174. MR378204

McNeil, A. J. (1997). Estimating the tails of loss severity distributions using extreme value theory. *ASTIN Bulletin* **27** 117–137.

Owen, A. B. (1988). Empirical likelihood ratio confidence intervals for a single functional. *Biometrika* **75** 237–249. MR946049

Owen, A. B. (1990). Empirical likelihood ratio confidence regions. *Ann. Statist.* **18** 90–120. MR1041387

Owen, A. B. (2001). *Empirical Likelihood*. Chapman and Hall, London.

Peng, L. (2001). Estimating the mean of a heavy tailed distribution. *Statist. Probab. Lett.* **52** 255–264. MR1838213

Qin, J. and Wong, A. (1996). Empirical likelihood in a semi-parametric model. *Scand. J. Statist.* **23** 209–219. MR1394654

Resnick, S. I. (1996). Discussion of the Danish data on large fire insurance losses. Preprint.

Resnick, S. I. (1997). Heavy tail modelling and teletraffic data. *Ann. Statist.* **25** 1805–1869. MR1474072

School of Mathematics
Georgia Institute of Technology
Atlanta, Georgia 30332-0160
USA
e-mail: peng@math.gatech.edu